\documentclass{amsart}
\usepackage{amssymb}
\makeatletter

\def\<#1\>{\left\langle#1\right\rangle}

\let\Re\undefined

\DeclareMathOperator{\Re}{Re}

\DeclareMathOperator{\Res}{Res}
\DeclareMathOperator{\Ind}{Ind}

\DeclareMathOperator{\ch}{ch}
\DeclareMathOperator{\id}{id}
\DeclareMathOperator{\ord}{ord}
\DeclareMathOperator{\vol}{vol}
\DeclareMathOperator{\meas}{meas}

\newcommand{\dd}[1]{\mathop{\mathrm{d}#1}}

\newcommand{\al}{\alpha}
\newcommand{\be}{\beta}
\newcommand{\ga}{\gamma}
\newcommand{\de}{\delta}
\newcommand{\ep}{\varepsilon}
\newcommand{\la}{\lambda}
\newcommand{\ze}{\zeta}
\newcommand{\wg}{\varpi}
\newcommand{\ph}{\varphi}
\newcommand{\om}{\omega}
\newcommand{\De}{\varDelta}
\newcommand{\Th}{\varTheta}
\renewcommand{\Pi}{\varPi}

\newcommand{\phs}{\ph^\circ}

\newcommand{\limp}{\mathrel{\enspace\Longrightarrow\enspace}}

\newcommand{\psm}[2][\Biggl(]
  {\left(\vphantom{#1}\begin{smallmatrix}#2\end{smallmatrix}\right)}

\def\abs{\@ifstar{\abs@star}{\abs@nostar}}
\newcommand{\abs@nostar}[1]{\lvert#1\rvert}
\newcommand{\abs@star}[1]{\left\lvert#1\right\rvert}
\def\Abs{\@ifstar{\Abs@star}{\Abs@nostar}}
\newcommand{\Abs@nostar}[1]{\lVert#1\rVert}
\newcommand{\Abs@star}[1]{\left\lVert#1\right\rVert}

\newcommand\set[1]{\Bigl\{#1\Bigr\}}

\newcommand{\cj}[1]{{\overline{#1}}}

\newcommand{\ade}{\mathbb{A}}
\newcommand{\ide}{\mathbb{J}}

\newcommand{\lqu}{\backslash}

\newcommand{\makegroup}[7][k]{%
  \newcommand{#3}[1][#1]{{#2_{##1}}}%
  \newcommand{#4}{{#2_\ade}}%
  \newcommand{#5}[1][#1]{{#3[##1]\lqu#4}}%
  \newcommand{#6}[1][v]{{#2_{##1}}}
  \newcommand{#7}{{#2_\infty}}}

\makegroup{G}{\Gk}{\Ga}{\Gq}{\Gv}{\Gin}
\makegroup{H}{\Hk}{\Ha}{\Hq}{\Hv}{\Hin}
\makegroup{M}{\Mk}{\Ma}{\Mq}{\Mv}{\Min}
\makegroup{N}{\Nk}{\Na}{\Nq}{\Nv}{\Nin}
\makegroup{P}{\Pk}{\Pa}{\Pq}{\Pv}{\Pin}
\makegroup{Q}{\Qk}{\Qa}{\Qq}{\Qv}{\Qin}
\makegroup{\Th}{\Sk}{\Sa}{\Sq}{\Sv}{\Sin}

\newcommand{\GL}{\mathrm{GL}}
\newcommand{\gO}{\mathrm{O}}

\newcommand{\SO}{\mathrm{SO}}

\let\frk\mathfrak

\newcommand{\CC}{\mathbb{C}}

\newcommand{\kx}{{k^\times}}

\newcommand{\oo}{\frk o}

\newcommand{\kext}{{\widetilde k/k}}
\newcommand{\Cx}{{\CC^\times}}

\renewcommand{\thmhead}[3]{\thmnumber{{\rm(#2)}\@ifnotempty{#1#3}%
    { }}\thmname{#1\@ifnotempty{#3}{\,---\,}}\thmnote{#3}}
\newcommand{\thref}{\eqref}

\newtheorem{prop}[equation]{Proposition}
\newtheorem{lem}[equation]{Lemma}
\newtheorem*{thm*}{Theorem}
\newtheorem*{prop*}{Proposition}
\newtheorem*{lem*}{Lemma}
\theoremstyle{definition}

\newtheorem*{defn*}{Definition}

\usepackage[lite]{amsrefs}

\begin{document}

\title[Periods of Eisenstein series of orthogonal groups of rank one]
  {Compact periods of Eisenstein series of orthogonal groups of rank one}
\author{Jo\~ao Pedro Boavida}
\address{Departamento de Matem\'atica\\
         Instituto Superior T\'ecnico\\Universidade T\'ecnica de Lisboa\\
         Av.\ Rovisco Pais\\1049--001 Lisboa, Portugal}
\email{joao.boavida@ist.utl.pt}

\begin{abstract}
Fix a number field $k$ with its adele ring $\ade$.
Let $G=\gO(n+3)$ be an orthogonal group of $k$-rank $1$ and $H=\gO(n+2)$ a
$k$-anisotropic subgroup.
We unwind the global period
\begin{equation*}
  (E_\ph,F)_H = \int_{\Hq}E_\ph\cdot\cj F
\end{equation*}
of a spherical Eisenstein series $E_\ph$ of $G$ against a cuspform $F$
of $H$ into an Euler product and evaluate the local factors at odd
primes.
\end{abstract}

\maketitle

\section*{Introduction}

Fix a number field $k$ with its adele ring $\ade$.  Equip $k^{n+3}$ with a
quadratic form $\<\;,\;\>$ with matrix
\begin{equation*}
  \begin{pmatrix}1&&\\&*&\\&&-1\end{pmatrix}
\end{equation*}
with respect to the orthogonal decomposition $k^{n+3}=(k\cdot
e_+)\oplus k^{n+1}\oplus(k\cdot e_-)$.  Let $G=\gO(n+3)$ and its
subgroups act always on the right.

In the same coordinates, write $k^{n+2}=(k\cdot e_+)\oplus k^{n+1}$
and let $H=\gO(n+2)$ be the fixer of $e_-$.  Let also $\De$ be the
discriminant of the restriction of $\<\;,\;\>$ to $k^{n+2}$ and let
$\chi$ be the corresponding quadratic character.

We consider only the case when $k^{n+2}$ is anisotropic; in
particular, $G$ has $k$-rank $1$ and $H$ is $k$-anisotropic.  With
$e=e_++e_-$, let $P$ be the rational parabolic stabilizing the
isotropic line $k\cdot e$.  Its unipotent radical $N^P$ is the fixer of
$e$ and $\{e_+,e_-\}^\perp=k^{n+1}$.  We choose the Levi component $M^P$ to 
be the stabilizer of $k^{n+1}$ and $\Th=\gO(n+1)$ to be its subgroup
fixing both $e_+$ and $e_-$.  Therefore, $M^P$ has the form
$\GL(1)\times\Th$.

We determine the period
\begin{equation*}
  (E_\ph,F)_H
  = \int_{\Hq} E_\ph\cdot \cj F
\end{equation*}
along $H$ of a spherical Eisenstein series $E_\ph$ on $G$ against a
cuspform $F$ of $H$.

Given a Hecke character $\om:\GL(1)\to\Cx$, extend it trivially to
$M^P$, and then, by right $N^P$--invariance, to $\ph_\om:\Pq\to\CC$.

\begin{prop*}
Let $E_\ph=E_{\om,\eta}$ be the Eisenstein series obtained from
$\ph_\om\cdot\eta$, where $\ph_\om$ is $\Th$-invariant and $\eta$ 
is a cuspform on $\Th$.  Then
\begin{equation*}
  (E_\ph,F)_H
  = \int_{\Sa\lqu\Ha} \ph_\om(h)\cdot (\eta,h\cdot F)_\Th\dd h.
\end{equation*}  
\end{prop*}
We shall see the inner period can be simplified further, essentially
reducing the original problem to the case where $\ph_\om=\ph_s$ is obtained
from a Hecke character with the same local parameter $s$ everywhere:
\begin{prop*}
Let $\ph=\ph_s\cdot\eta$ be as above and $\al=(n+1)s$.  Up to
correction factors at finitely many places, the period against the
constant function is
\begin{align*}
  (E_\ph,1)_H
  &=\frac{(\eta,1)_\Th \cdot \vol(\Th\lqu\Th K^H) \cdot \ze_k(\al-n)}
      {L_k(\al-\lfloor\frac n2\rfloor,\chi)},
    &&\text{if $n$ is odd;}\\
  (E_\ph,1)_H
  &=\frac{(\eta,1)_\Th \cdot \vol(\Th\lqu\Th K^H) \cdot \ze_k(\al-n)}
      {\ze_k(2\al-n) / L_k(\al-\frac n2,\chi)},
    &&\text{if $n$ is even.}
\end{align*}
The set of bad primes consists of all archimedean primes and all factors
of $2\De$.
\end{prop*}

Clearly, this period is nonzero if and only if $(\eta,1)_\Th\ne0$.
Section \ref{s:gen} discusses $(E_\ph,F)_H$ in the general case.

These findings fit into the Gross--Prasad conjecture on periods of
$\SO(n+1)$--automorphic forms over $\SO(n)$.  In the case at hand, the
conjecture \cites{GrPr1, GrPr2, GrPr3} predicts that a representation
of $\gO(n+2)$ occurs in a representation of $\gO(n+3)$ if and only if the
corresponding tensor product $L$-function is nonzero on $\Re s=\frac12$.
Ichino and Ikeda \cite{IcIk} discuss further details and broader
context is provided in papers by Jacquet, Lapid, Offen, and/or
Rogawski \cites{JaLaRo,LaRo,LaOf}, and surveyed by Gross, Reeder
\cite{GrRe}, and Jiang \cite{Ji}.

These same periods (called there global Shintani functions) were used
by Murase and Sugano \cite{MuSu} to obtain and study integral expressions
for standard $L$-functions of the orthogonal group.

In section~\ref{s:setup}, we unwind the period into an Euler product, 
obtaining
\begin{equation*}
  (E_\ph,F)_H
  = \int_{\Sk\lqu\Ha}\ph\cdot \cj F
  = \int_{\Sa\lqu\Ha}\ph_\om(h)\;
      \int_{\Sq}\eta(\theta)\cdot\cj F(\theta h)\dd\theta\dd h.
\end{equation*}
For several sections, we focus on $\eta=1$ and $F=1$, in which case
\begin{equation*}
  (E_\ph,1)_H
  = \vol(\Sq)\cdot\prod_v\int_{\Sv\lqu\Hv}\ph_{\om,v}.
\end{equation*}

Because $\ph$ is spherical, the local integral is trivial at
anisotropic places.  At isotropic places, the local evaluation entails
a choice of local parametrization for $\Sv\lqu\Hv$, as well as a
normalization of the measure.  In section~\ref{s:odd}, we consider the
specifics at odd primes.  Notably, at good primes it is possible to
reduce the evaluation to $n=1$ or $n=2$.  We articulate the full
details in the subsequent sections.  In section \ref{s:odd:bad}, we
determine the correction factors at bad odd primes.  In
section~\ref{s:gen}, we resume the general case and conclude that
$(E_\ph,F)_H$ is obtained from $(\eta,F)_\Th$ as well as the local
integrals (obtained in the preceding sections) of $\ph_\om$ against a
degenerate principal series of $H$.  We show also, following Garrett
\cite{Ga}, that the periods of a cuspform are nonzero only if, almost
everywhere, the cuspform generates locally the very same sort of
degenerate principal series generated by the Eisenstein series we
consider.  Finally, in the appendix, we use the same methods
to evaluate (the local factors of) the constant term of $E_\ph$.

Throughout, we rely on
$\gO(n+1)\supset\gO(n)$ being a Gelfand pair, a fact conjectured by
J.~Bernstein and established independently by Kato, Murase, and Sugano
\cite{KaMuSu} and Aizenbud, Gourevitch, and Sayag \cite{AiGoSa}.

\subsection*{Acknowledgements}

This paper is based on part of the author's doctoral dissertation, done 
under the supervision of Paul Garrett.  It is influenced by discussions with and talks by him.  The author wishes also to acknowledge the helpful advice 
of a referee.

\tableofcontents

\section{Setup}\label{s:setup}

Recall we fixed a number field $k$ with adele ring $\ade$, 
and a quadratic form $\<\;,\;\>$ with matrix
\begin{equation*}
  \begin{pmatrix}1&&\\&*&\\&&-1\end{pmatrix}
\end{equation*}
with respect to the decomposition $k^{n+3}=(k\cdot e_+)\oplus
k^{n+2}\oplus(k\cdot e_-)$.  We set $e=e_++e_-$ and named the
following groups of isometries:
\begin{align*}
  G&=\gO(n+3), &&\text{the isometry group of }
    \begin{pmatrix}*&*&*\end{pmatrix};\\
  H&=\gO(n+2), &&\text{the isometry group of }
    \begin{pmatrix}*&*&\phantom{*}\end{pmatrix};\\
  \Th&=\gO(n+1), &&\text{the isometry group of }
    \begin{pmatrix}\phantom{*}&*&\phantom{*}\end{pmatrix};\\
  P&\subset G,&&\text{the $k$-parabolic stabilizing }k\cdot e;\\
  N^P&\subset P,&&\text{its unipotent radical, fixing $e$ and }
    \begin{pmatrix}\phantom{0}&*&\phantom{0}\end{pmatrix};\\
  M^P&\subset P,&&\text{the Levi component stabilizing }
    \begin{pmatrix}\phantom{0}&*&\phantom{0}\end{pmatrix}. 
\end{align*}
The modular function of $P$ is given by $\de_P(p)=\abs t^{n+1}$ when 
$e\cdot p=e/t$.  In particular, $\de_P^s(p)=\abs t^\al$, with 
$\al=(n+1)s$.

Recall $M^P$ has the form $\GL(1)\times\Th$.  Given a Hecke
character $\om:\kx\lqu\ide\to\Cx$, extend it trivially to $M^P$, and
then, by right $N^P$--invariance, to $\ph_\omega:\Pq\to\CC$.  Given also
a cuspform $\eta:\Sq\to\CC$ normalized by $\eta(1)=1$ and spherical, define
$\ph:\Pq\to\CC$ by $\ph(\theta g)=\eta(\theta)\cdot\ph_\omega(g)$,
where $\theta\in\Th$ and $g\in\GL(1)\cdot N^P$.

Let
\begin{equation*}
  E_\ph(g)
  = E_{\om,\eta}(g)
  = \sum_{\ga\in\Pk\lqu\Gk}\ph(\ga g).
\end{equation*}
be the \emph{spherical Eisenstein series} of $G$ associated with
$\ph$, and let $F$ be a cuspform of $H$.
We evaluate the period 
\begin{equation*}
  (E_\ph,F)_H
  = \int_{\Hq} E_\ph\cdot\cj F
  = \int_{\Hq} \sum_{\ga\in\Pk\lqu\Gk}\ph(\ga h) \cdot\cj F(h)\dd h.
\end{equation*}

By Witt's lemma, $\Pk\lqu\Gk$ is the space of isotropic lines in 
$k^{n+3}$ and $\Hk$ acts transitively on that space (here we 
rely on the $k$--Witt index being $1$, so that no 
$k$--isotropic lines are orthogonal to $e_-$).  
Because $\Th=H\cap P=\gO(n+1)$ is the fixer of $e_+$ and $e_-$,
\begin{equation*}
  (E_\ph,F)_H
  = \int_{\Hq}\sum_{\ga\in\Sk\lqu\Hk}\ph(\ga h)\cdot\cj F(h)\dd h
  = \int_{\Sk\lqu\Ha}\ph\cdot\cj F.
\end{equation*}
As $\ph(\theta g)=\eta(\theta)\cdot\ph_\om(\theta g)$ and $\ph_\om$ 
is left $\Th$--invariant,
\begin{equation*}
  (E_\ph,F)_H
  = \int_{\Sa\lqu\Ha} \ph_\om(h)
      \int_{\Sq}\eta(\theta)\cdot\cj F(\theta h)\dd\theta\dd h.
\end{equation*}
Because $H\supset\Th$ is a Gelfand pair \cite{AiGoSa} and $\eta$ is
spherical, the inner integral can be expressed as
\begin{equation*}
  \int_{\Sq}\eta(\theta)\cdot\cj F(\theta h)\dd\theta
  = (\eta,F)_\Th\cdot f(h),
\end{equation*}
where $f$ is a spherical vector of
$\Ind_\Th^H1$ normalized by $f(1)=1$.  Therefore,
\begin{equation}\label{e:auto-H}
  (E_\ph,F)_H
  = (\eta,F)_\Th \cdot \int_{\Sa\lqu\Ha}\ph_\om\cdot f.
\end{equation}

From now on and until section \ref{s:gen}, we restrict our attention
to $\eta=1$ and $F=1$ (thus, also $\ph=\ph_\om$ and $f=1$) and work
locally.  Therefore, we omit mention of the place $v$ whenever
possible.

\subsection*{Anisotropic places}

At anisotropic places, the local integral is $\vol(\Th\lqu H)$.  
(Recall we are dealing with spherical functions only.)

\subsection*{Isotropic places}

Choose a hyperbolic pair $x$, $x'$ in $k^{n+2}$ so that 
$e_+\in k\cdot\{x,x'\}$ and change coordinates so that
the restricted quadratic form has matrix 
\begin{equation*}
  \begin{pmatrix}&&1\\&B&\\1&&\end{pmatrix}
\end{equation*}
with respect to the orthogonal decomposition 
$k^{n+2}=(k\cdot x')\oplus k^n\oplus(k\cdot x)$.  Recall that 
$H=\gO(n+2)$ and that 
$\Th=\gO(n+1)$ is the subgroup of $H$ fixing $e_+$.  We are to evaluate
\begin{equation*}
  \int_{\Th\lqu H}\ph.
\end{equation*}

Let $Q\subset H$ be the parabolic stabilizing the line $k\cdot x$;
we claim that $\Th Q$ is the open orbit of $\Th\lqu H/Q$.  
By Witt's lemma, $\Th\lqu H$ can be identified with the homogeneous 
space of vectors $y$ with $\<y,y\>=1$ 
and the orbit of $e_+$ under the action of $Q$ includes all such $y$ 
not orthogonal to $x$:
the only non-trivial requirement is that $\<yq,xq\>=\<y,x\>$, which 
can be achieved by choosing $xq=\mu x$ and 
$\<y,x\>=\<yq,xq\>=\mu\<yq,x\>$.

Therefore, 
\begin{equation}\label{e:quo}
  \int_{\Th\lqu H} \text{function}(h)\dd h
  = \int_{\Th\lqu\Th Q} \text{function}(q)\dd q
  = \int_{(\Th\cap Q)\lqu Q} \text{function}(q)\dd q.
\end{equation}
Here, $\Th\cap Q=\gO(n)$ is the fixer of $x$ and $x'$.  Set
\begin{equation*}
  m_\la = \begin{pmatrix}\la&&\\&\id&\\&&1/\la\end{pmatrix}
\qquad\text{and}\qquad
  n_a = \begin{pmatrix}1&a&-\frac12B(a)\\&\id&*\\
     \phantom{\frac12B(a)}&&1\end{pmatrix}.
\end{equation*}
With $M_* = \bigl\{m_\la\bigr\}$, we have
$\Th\cap Q=\Bigl\{\psm[.]{1&&\\&*&\\&&1}\Bigr\}$, $N^Q=\bigl\{n_a\bigr\}$, 
$M^Q=(\Th\cap Q)\cdot M_*$, and $Q=M^Q\cdot N^Q$.
The elements of $(\Th\cap Q)\lqu Q$ can be expressed as $n_a\cdot m_\la$
and $\de_Q(m_\la)=\abs\la^n$.  Moreover,
\begin{equation*}
  \dd{(n_a\cdot m_\la)}
  = \dd a\dd\la
\end{equation*}
(with $\dd\la$ multiplicative and $\dd a$ additive) is a 
\emph{right}-invariant measure.

From this point onward, we must consider each place
separately.

\section{Odd primes}\label{s:odd}

Fix (and suppress) an odd place $v$.  Choose coordinates so that the
maximal compact open subgroup 
$K$ of $H$ consists of the orthogonal matrices with integral
entries.  With $\Phi$ being the characteristic function of $\oo^{n+3}$
and $\om$ a character of $\kx$, define, for this section only,
\begin{equation*}
  \ph(g)
  = \int_{\kx}\om(t)\cdot\Phi(te\cdot g)\dd t
  \qquad\text{and}\qquad
  \phs
  = \frac{\ph}{\ph(1)}.
\end{equation*}
Note that this is a change in notation: what was denoted $\ph_\om$ in
the last section, will be denoted $\phs$ in this one.  We address the
case $\om=\abs{\,\cdot\,}^\al$.

We set notation for the next several sections.  Recall that $\al=(n+1)s$.
The discriminant enters the picture as $\ep=\chi(\De)$, where
$\chi$ is the non-trivial character on $\oo^\times/{\oo^\times}^2$
and $\widetilde k=k[x]/\<x^2-\De\>$ is the quadratic extension.
At bad primes (those dividing the discriminant), we set $\ep=0$.
We write
\begin{equation*}
  Z(\al) = \int_{\kx\cap\oo}\abs t^\al\dd t
  = \frac1{1-q^{-\al}}
  \qquad\text{and}\qquad
  L(\al,\chi) = \frac1{1-\ep q^{-\al}}.
\end{equation*}
We write also
\begin{equation*}
  Z_\kext(\al)=Z(\al)\; L(\al,\chi).
\end{equation*}
We set $\abs t=q^{-T}$ and $a=q^{-\al}$ throughout.  We shall use
\begin{equation*}
  \int_{\kx\cap\oo}|t|^\al \; u^T\dd t
  = \int_{\kx\cap\oo}(au)^T\dd t
  = \sum_{T\ge0}(au)^T = \frac1{1-au} = Z(\al-\log_q u).
\end{equation*}

In this section, we express the \emph{local integral}
\begin{equation*}
  \int_{\Th\lqu H}\phs
\end{equation*}
in terms of the functions $X$ and $\Pi$ that we are about to define.  
We prove also it is sufficient to evaluate those functions (and
the corresponding local integrals) for small $n$.

\begin{defn*}
With $z=q^{-\be}$, we define
\begin{align*}
  X_\ell^B(\rho)
  &=\meas\set{a\in\oo^n:\frac{B(a)-\rho}2=0\bmod\wg^\ell};\\[2\jot]
  X^B(\beta;\rho)
  &=\sum_{\ell\ge0}z^\ell X_\ell^B(\rho);\qquad\text{and}\\
  \Pi^B(\al,\be)
  &= \int_{\kx\cap\oo}\abs t^\al \; X^B(\be;t^2)\dd t.
\end{align*}
(When there is no risk of ambiguity, we suppress $B$ or $\rho$, 
or use $n$ instead of $B$.)
\end{defn*}

\begin{prop}\label{p:odd:lpi}
The local integral is obtained from
\begin{equation*}
  \int_{\Th\lqu H}\phs
  = \frac{\vol(\Th\lqu\Th K)\cdot\Pi^n(\al-\be-n,\be)}{Z(\al)\cdot X^n(0;1)}
\end{equation*}
with $\be=0$.
\end{prop}

\begin{proof}
To account for the normalization implied in the integral \eqref{e:quo}, 
we consider also the integral of $\Phi(e\cdot h)$, which
happens to be, at odd primes, the characteristic function of 
$\Th\lqu\Th K$.  (Indeed, for $h\in H$, both or neither of $e\cdot h$ and
$e'\cdot h$ are in $\oo^{n+3}$.  The statement follows at once from the
Cartan decomposition in $G$.)

That allows us to express the local integral (still 
without $\abs\la^\be$) as
\begin{equation}\label{e:odd:lp}
  \vol(\Th\lqu\Th K)\cdot
    \frac{\int_{\Th\lqu H}\phs}{\int_{\Th\lqu H}\ch_{\Th\lqu\Th K}}
  = \frac{\vol(\Th\lqu\Th K)\cdot\int_{(\Th\cap Q)\lqu Q} \ph}
      {\ph(1)\cdot\int_{(\Th\cap Q)\lqu Q} \ch_{\Th\lqu\Th K}}.
\end{equation}

With a new parameter $\be=0$, the integral in 
the numerator of \eqref{e:odd:lp} is
\begin{equation}\label{e:odd:triple}
  \int_{(\Th\cap Q)\lqu Q}\ph\cdot\abs\la^\be
  = \int_{\kx}\int_{\kx}\int_{k^n}
      \abs t^\al \; \abs\la^\be \;
      \Phi(te\cdot n_a\cdot m_\la)\dd a\dd\la\dd t.
\end{equation}

Recall that we chose coordinates so that $e_+$ be a linear combination
of $x$ and $x'$.  We specify further that $e_+=x'+\frac12x$.  
Noting that
\begin{equation*}
  e\cdot n_a\cdot m_\la
  = (e_+ + e_-)\cdot n_a\cdot m_\la
  = e_+\cdot n_a\cdot m_\la + e_-,
\end{equation*}
we have (in $k^{n+2}$)
\begin{equation*}
  e_+\cdot n_a\cdot m_\la
  = \begin{pmatrix}1&0&\tfrac12\end{pmatrix}\cdot n_a\cdot m_\la
  = \begin{pmatrix}\la&a&\tfrac1{2\la}\bigl(1-B(a)\bigr)\end{pmatrix}
\end{equation*}
and (in $k^{n+3}$)
\begin{equation*}
  te\cdot n_a\cdot m_\la
  = \Bigl(\la t,at,\tfrac1{2\la t}\bigl(t^2-B(at)\bigr),t\Bigr).
\end{equation*}

Therefore, continuing \eqref{e:odd:triple} and after a change of variables,
\begin{equation*}\begin{split}
  \int_{(\Th\cap Q)\lqu Q}\ph \cdot \abs\la^\be
  &=\int_{\kx}\int_{\kx}\int_{k^n}
      \abs t^{\al-\be-n} \; \abs\la^\be \; \Phi
      \Bigl(\la,a,\tfrac1{2\la}\bigl(t^2-B(a)\bigr),t\Bigr)
      \dd a\dd\la\dd t\\
  &=\int_{\kx\cap\oo}\abs t^{\al-\be-n}
      \int_{\kx\cap\oo}\abs\la^\be
      \int_{\oo^n}\ch_\oo \Bigl(\frac{t^2-B(a)}{2\la}\Bigr)
      \dd a\dd\la\dd t.
\end{split}\end{equation*}

Consider now the integral in the denominator of \eqref{e:odd:lp}.
We have
\begin{equation*}\begin{split}
  \int_{(\Th\cap Q)\lqu Q}\ch_{\Th\lqu\Th K}
  &=\int_{\kx}\int_{k^n}\Phi(e\cdot n_a\cdot m_\la)\dd a\dd\la\\
  &=\int_{\kx}\int_{k^n}
    \Phi\Bigl(\la,a,\tfrac1{2\la}\bigl(1-B(a)\bigr),1\Bigr)\dd a\dd\la
  = X(0;1).\qedhere
\end{split}\end{equation*}
\end{proof}

\subsection*{Dimension reduction}

By taking hyperbolic planes away, we can simplify the evaluation of
\thref{p:odd:lpi} significantly.  The argument is valid at even primes too.

We write vectors in $k^{n'}=k^{n+2}$ as $a'=(x,a,y)$ and the quadratic form
as $B'(a')=B(a)-2xy$, where $a\in k^n$ and $B(a)$ is the restriction to 
$k^n$.  Note that $B'$ and $B$ have the same discriminant.
Also, we abbreviate $X^{B'}_\ell(\rho)=X'_\ell(\rho)$, 
$X^{B'}(\be;\rho)=X'(\be;\rho)$, $X^B_\ell(\rho)=X_\ell(\rho)$, and
$X^B(\be;\rho)=X(\be;\rho)$.  In particular,
\begin{align*}
  X'_\ell
  &=\meas\set{(x,a,y)\in\oo^{n+2}:\frac{B(a)-\rho}2-xy=0\bmod\wg^\ell};\\
  X_\ell
  &=\meas\set{a\in\oo^n:\frac{B(a)-\rho}2=0\bmod\wg^\ell}.
\end{align*}

\begin{prop}\label{p:odd:reduce}
With $z=q^{-\be}$, we have
\begin{equation*}
  X'(\be;\rho)
  = \frac{Z(\be+1)}{Z(\be+2)}\cdot X(\be+1;\rho).
\end{equation*}
Moreover, if there is a hyperbolic subspace with dimension $2k$ and $n=m+2k$,
\begin{align*}
  X^{m+2k}(\be;\rho)
  &=\frac{Z(\be+1)}{Z(\be+k+1)}\cdot X^m(\be+k;\rho),\\
  \Pi^{m+2k}(\al,\be)
  &=\frac{Z(\be+1)}{Z(\be+k+1)}\cdot\Pi^m(\al,\be+k),
\end{align*}
and the local integral is
\begin{equation*}
  \frac{\vol(\Th\lqu\Th K)\cdot\Pi^m(\al-\be-n,\be+k)}{Z(\al)\cdot X^m(k;1)}.
\end{equation*}
\end{prop}

\begin{proof}
If $\rho^{(2)}=(B(a)-\rho)/2\bmod\wg^\ell$, then the set
\begin{equation*}
  \set{a\in\oo^n : \frac{B(a)-\rho}2-xy=0\bmod\wg^\ell}
\end{equation*}
breaks into
\begin{equation*}
  \bigcup_{N\ge0}
    \set{a\in\oo^n : \ord\rho^{(2)}=N ; \rho^{(2)}-xy=0\bmod\wg^\ell}.
\end{equation*}

We shall prove in proposition \thref{p:odd:X2} that, for
an isotropic plane, $X^2_\ell(2\rho^{(2)})$ depends only on
$\ord\rho^{(2)}$ (we express it thus so that the argument be valid
at even primes too). That allows us to write
\begin{equation*}
  X'_\ell(\rho)
  = \sum_{N\ge0}\bigl(X_N(\rho)-X_{N+1}(\rho)\bigr)\cdot X^2_\ell(2\wg^N)
\end{equation*}
and
\begin{equation*}
  X'(\be;\rho)
  = \sum_{\ell\ge0}z^\ell X'_\ell(\rho)
  = \sum_{N\ge0}(X_N-X_{N+1})\cdot X^2(2\wg^N).
\end{equation*}
Relying on proposition \thref{p:odd:X2} (which, for the isotropic case
and even primes \emph{is} about $X^2(2\wg^N)$)
and setting $w=zq^{-1}$, we obtain
\begin{equation*}
  X'(\be;\rho)
  = \frac{Z(\be+1)}{Z(\be+2)}\cdot X(\be+1;\rho).\qedhere
\end{equation*}
\end{proof}

\subsection*{Strategy overview}

The main task left is to determine $X(\be;t^2)$ and 
$\Pi(\al,\be)$ for anisotropic forms.  It will prove convenient to use
$z=q^{-\be}$ and $a=q^{-\al}$ consistently.  We will also use $w=zq^{-1}$, as
well as $u=q^{-2\be-m}$ (with $m$ being the dimension 
of the anisotropic space under consideration), so that when $\be$ is set to 
the Witt index (of $k^n$) we obtain $u=q^{-n}$.

According to propositions \thref{p:odd:lpi} and \thref{p:odd:reduce}, 
if $n=m+2k$ (where $k$ denotes now the Witt index) the local integral is
\begin{multline*}
  \frac{\vol(\Th\lqu\Th K)\cdot\Pi^n(\al-\be-n,\be)}
    {Z(\al)\cdot X^n(0; 1)}\\
  = \frac{\vol(\Th\lqu\Th K)\; Z(\be+1) \; Z(k+1) \; \Pi^m(\al-\be-n,\be+k)}
    {Z(\al) \; Z(\be+k+1) \; Z(1) \; X^m(k; 1)}.
\end{multline*}
In particular, for $\be=0$ (which is what we need for the period), we obtain
\begin{equation*}
  \frac{\vol(\Th\lqu\Th K)\cdot \Pi^m(\al-n, k)}{Z(\al)\cdot X^m(k; 1)}.
\end{equation*}

We shall find that $X^m(\be;t^2)$ is a rational function of $z$ (which we will
often write also with $w$ or $u$), whose dependence on $t^2$ is in fact 
a dependence on $T=\ord t$ only.  Moreover,
\begin{equation*}
  \Pi^m(\al, \be)
  = \int_{\kx\cap\oo} \abs t^\al \; X^m(\be; t^2)\dd t
  = \sum_{T\ge0} a^T X^m(\be; \wg^{2T}).
\end{equation*}
This will be a rational function $R(z,w,u,a)$ of $z$ (or $w$ or $u$) and $a$, 
whose factors all have the form $1\pm q^*$ after $z$, $w$, $u$ and $a$ are 
substituted (and $a$ occurs only with non-negative exponents).  At $a=0$, we 
obtain the coefficient of $T=0$.  That is, $X^m(\be; 1)=R(z,w,u,0)$ 
consists of the factors of $R$ that do not involve $a$.

In other words, if $R^*(z,w,u,a)$ consists of the 
factors of $R(z,w,u,a)$ that \emph{do} involve $a$, then 
\begin{equation*}
  \frac{\Pi^m(\al,\be)}{X^m(\be;1)}
  = R^*(z,w,u,a).
\end{equation*}
That is, the local period is
\begin{equation}\label{e:odd:master}
  \vol(\Th\lqu\Th K)\cdot
    \frac{R^*\!\left(q^{-k},q^{-k-1},q^{-n},q^{-(\al-n)}\right)}{Z(\al)}.
\end{equation}
(We will often suppress $\vol(\Th\lqu\Th K)$ from this point onward.)

In the prototypical scenario, we might obtain an expression like
\begin{equation}\label{e:proto:X}
  X(\be;t^2)
  = \frac{\text{fn}(z,w,u)}{1-u} \;
    \Biggl[1 - \frac{u-v}{1-v} u^T\Biggr],
\end{equation}
(where $v$ might be $w$, $\ep wq^{-1}$, $-z$, or a similar term).  Then
\begin{equation*}
  \Pi(\al,\be)
  = \frac{\text{fn}(z,w,u)}{1-u} \;
    \Biggl[\frac1{1-a} - \frac{u-v}{(1-v)(1-au)}\Biggr]
  = \frac{\text{fn}(z,w,u) \; (1-av)}{(1-v)(1-a)(1-au)}.
\end{equation*}
In light of what we said above, the local period \eqref{e:odd:master} would 
then be
\begin{equation}\label{e:proto:per}
  \frac1{Z(\al)}\cdot\frac{1-av}{(1-a)(1-au)}
  = \frac{Z(\al-n)\;Z(\al)}{Z(\al)}\cdot(1-av)
  = \frac{Z(\al-n)}{1/(1-av)}.
\end{equation}

\section{Good odd primes---$n=2k+1$}

We consider first $n=1$.
The quadratic form on $k^1$ can be expressed as $\De x^2$.

\begin{prop}\label{p:odd:X1}
With $z=q^{-\be}$, $w=zq^{-1}$, $u=z^2q^{-1}$, and $\De$ a unit, we have
\begin{equation*}
  X^1(\be;t^2)
  = \frac{1+w}{1-u} \; \Biggl[1 - \frac{u-\ep w}{1-\ep w} u^T\Biggr].
\end{equation*}
\end{prop}

\begin{proof}
We first evaluate
\begin{equation*}
  X_\ell(t^2)
  = \meas\set{x\in\oo^1 : \De x^2 = t^2\bmod\wg^\ell}.
\end{equation*}

If $t^2=0\bmod\wg^\ell$ (i.e.,\ if $\ell\le2T$), we have
\begin{equation*}
  X_\ell = X_\ell(0)
  = \meas\set{x\in\oo : x^2=0\bmod\wg^\ell}
  = q^{-\lceil\ell/2\rceil}
\end{equation*}
and
\begin{equation}\label{e:zlqfrac}
  \sum_{0\le\ell\le2T} z^\ell q^{-\lceil\ell/2\rceil}
  = 1 + \sum_{0<k\le T} (z^{-1}+1) z^{2k} q^{-k}
  = \frac{1+w}{1-u} - \frac{w+u}{1-u} u^T.
\end{equation}

On the other hand, if $t^2\ne0\bmod\wg^\ell$ (i.e.,\ if $\ell>2T$), 
we have two distinct situations: either $\De x^2-t^2$ is isotropic ($\ep=1$) 
or it isn't ($\ep=-1$).

When $\ep=-1$, the form $\De x^2-t^2$ is anisotropic.
In that case, we show in lemma \thref{l:odd:aniso} below
that $\De x^2-t^2=0\bmod\wg^\ell$ requires $t^2=0\bmod\wg^\ell$.
Therefore, $X_\ell=0$ for $t^2\ne0\bmod\wg^\ell$.

When $\ep=1$, we have
\begin{equation*}
  \De x^2-t^2=0\bmod\wg^\ell
  \limp
  x^2=t^2\bmod\wg^\ell
  \limp
  \abs x = \abs t = q^{-T},
\end{equation*}
leading to $(x-t)(x+t)=0\bmod\wg^\ell$ and $x=\pm t\bmod\wg^{\ell-T}$.
That is,
\begin{equation*}
  X_\ell = \meas\set{x\in\oo:x=\pm t\bmod\wg^{\ell-T}}
  = 2q^{-\ell+T}
\end{equation*}
and
\begin{equation*}
  \sum_{\ell>2T}z^\ell X_\ell
  = \sum_{\ell>2T} z^\ell 2q^{-\ell+T}
  = \frac{2z^{2T+1}q^{-T-1}}{1-zq^{-1}}
  = \frac{2zq^{-1}(z^2q^{-1})^T}{1-zq^{-1}}
  = \frac{(1+\ep)wu^T}{1-w}.
\end{equation*}

Taking this last expression as well as \eqref{e:zlqfrac} into account, 
and combining the cases $\ep=\pm1$, we 
obtain the result stated.
\end{proof}

\begin{lem}\label{l:odd:aniso}
If $k$ is the local field at an odd prime and 
the form $B(x)=\sum a_ix_i^2$ is anisotropic, then
\begin{equation*}
  \abs{B(x)} = \max\set{\abs*{a_i x_i^2}}.
\end{equation*}
\end{lem}

\begin{proof}
If some coordinate had 
$\abs*{a_ix_i^2}>\abs{B(x)}$, then Hensel's lemma would yield a new vector 
$x'$ with $\ord x_i=\ord x'_i$ and $B(x')=0$.
\end{proof}

\begin{prop}\label{p:odd:per-even}
If $\De$ is a unit and $n=2k+1$, the local period is
\begin{equation*}
  \vol(\Th\lqu\Th K)\cdot\frac{Z(\al-n)}{L(\al-k,\chi)}
  = \vol(\Th\lqu\Th K)\cdot\frac{Z\bigl((n+1)(s-1)+1\bigr)}
      {L\bigl((n+1)(s-\tfrac12)+1,\chi\bigr)}.
\end{equation*}
\end{prop}

\begin{proof}
Let $z=q^{-\be}$, $w=zq^{-1}$, $u=z^2q^{-1}$ and $a=q^{-\al}$.  In 
proposition \thref{p:odd:X1}, we saw $X^1$ fits the prototype
\eqref{e:proto:X} with $v=\ep w$.  According to \eqref{e:proto:per}, the 
local period is
\begin{equation*}
  \frac{Z(\al-n)}{1/(1-\ep aw)}
  = \frac{Z(\al-n)}{L(\al-n+k+1,\chi)}.
\end{equation*}
The result follows from $\al=(n+1)s$.
\end{proof}

\section{Good odd primes---$n=2k+2$}

We consider $n=2$ first.

Expressing the form in terms of $x$ and $y$, we may assume $B=x^2-\De y^2$ 
and note that $X_\ell^B(\rho)$ depends only on
$\ord\rho$.  Indeed, any unit can be expressed in the form $\eta=a^2-\De b^2$ 
(this is a consequence of Hensel's lemma).  But
\begin{equation*}
  (ax+\De by)^2-\De(bx+ay)^2 = (a^2-\De b^2)(x^2-\De y^2) = \eta\,(x^2-\De y^2)
\end{equation*}
and the matrix $\bigl(\begin{smallmatrix}a&\De b\\b&a\end{smallmatrix}\bigr)$
is invertible in $\oo^2$ if its determinant $a^2-\De b^2$ is a unit.

\begin{prop}\label{p:odd:X2}
With $z=q^{-\be}$, $w=zq^{-1}$, $u=z^2q^{-2}$, and $\De$ a unit, we have
\begin{equation*}
  X^2(\be;t^2)
  = \frac{(1-\ep wq^{-1})(1-\ep wu^T)}{(1-w)(1-\ep w)}.
\end{equation*}
\end{prop}

\begin{proof}[Proof of proposition \thref{p:odd:X2} in anisotropic case]
This is the case $\ep=-1$.  We assumed the form is $x^2-\De y^2$.
We evaluate
\begin{equation*}
  X_\ell(t^2)
  = \meas\set{(x,y)\in\oo^2:x^2-\De y^2=t^2\bmod\wg^\ell},
\end{equation*}
for $\abs t=q^{-T}$.

When $t^2=0\bmod\wg^\ell$ (i.e.,\ $\ell\le2T$), the anisotropy lemma
\thref{l:odd:aniso} leads to
\begin{equation*}
  X_\ell(t^2)
  = \meas\set{(x,y)\in\oo^2:x^2=y^2=0\bmod\wg^\ell}
  = \bigl(q^{-\lceil\ell/2\rceil}\bigr)^2=q^{-2\lceil\ell/2\rceil}.
\end{equation*}
We saw in \eqref{e:zlqfrac} 
(using $q^{-2}$, $wq^{-1}=zq^{-2}$, and $u=z^2q^{-2}$ \emph{here} 
in place of $q^{-1}$, $w=zq^{-1}$, and $u=z^2q^{-1}$ \emph{there}) that
\begin{equation*}
  \sum_{0\le\ell\le2T}z^\ell X_\ell
  = \sum_{0\le\ell\le2T}z^\ell q^{-2\lceil\ell/2\rceil}
  = \frac{1+wq^{-1}}{1-u} - \frac{wq^{-1}+u}{1-u} u^T.
\end{equation*}

When $t^2\ne0\bmod\wg^\ell$ (i.e.,\ $\ell>2T$), again by 
anisotropy,
\begin{equation*}\begin{split}
  X_\ell(t^2)
  &=\meas\set{(x,y)\in\oo^2:x^2-\De y^2=t^2\bmod\wg^\ell}
  =\abs t^2\cdot X_{\ell-2T}(1).
\end{split}\end{equation*}

For each fixed $\ell$,
\begin{equation*}\begin{split}
  \sum_{\text{unit $\eta\bmod\wg^\ell$}}X_\ell(\eta)
  &=\sum_{\text{unit $\eta\bmod\wg^\ell$}}
    \meas\set{(x,y)\in\oo^2:x^2-\De y^2=\eta\bmod\wg^\ell}\\
  &=\meas\set{(x,y)\in\oo^2:x^2-\De y^2=\text{unit}}
  = 1-q^{-2}.
\end{split}\end{equation*}
(We relied on anisotropy for the last step.)  Because 
$X_\ell(\eta)=X_\ell(1)$, we see that
\begin{equation*}
  X_\ell(1)
  = \frac{1-q^{-2}}{\#\{\text{units $\eta\bmod\wg^\ell$}\}}
  = (1+q^{-1})q^{-\ell}
\end{equation*}
and also
\begin{equation*}
  X_\ell(t^2)
  = q^{-2T}\cdot X_{\ell-2T}(1)
  = (1+q^{-1})q^{-\ell}.
\end{equation*}
This leads to
\begin{equation*}
  \sum_{\ell>2T}z^\ell X_\ell
  = \sum_{\ell>2T} w^\ell(1+q^{-1})
  = \frac{w+wq^{-1}}{1-w}u^T.
\end{equation*}

Combining the two subsums (over $\ell\le2T$ and over $\ell>2T$) we obtain
the stated conclusion for the anisotropic case.
\end{proof}

\begin{proof}[Proof of proposition \thref{p:odd:X2} in isotropic case]
This is the case $\ep=1$.  We can assume the form is $xy$.
We will evaluate 
\begin{equation*}
  X_\ell(\rho)
  = \meas\set{(x,y)\in\oo^2: xy=\rho\bmod\wg^\ell}.
\end{equation*}
We set $N=\ord\rho$.

If $\rho=0\bmod\wg^\ell$ (i.e.,\ $\ell\le N$), we have
\begin{equation*}\begin{split}
  X_\ell
  &=\meas\set{x=0\bmod\wg^\ell}
    +\sum_{0\le k<\ell}
      \meas\set{\ord x=k; y=0\bmod\wg^{\ell-k}}\\
  &=q^{-\ell}+\sum_{0\le k<\ell}q^{-k}(1-q^{-1})q^{-\ell+k}
  = q^{-\ell}+\ell q^{-\ell}(1-q^{-1})
\end{split}\end{equation*}
and
\begin{equation*}
  \sum_{0\le\ell\le N}z^\ell X_\ell
  = \sum_{0\le\ell\le N}\Bigl(w^\ell+\ell w^\ell(1-q^{-1})\Bigr)
  = \frac{1-w^{N+1}}{1-w}
    + (1-q^{-1})\sum_{0\le\ell\le N}\ell w^\ell.
\end{equation*}
We observe in passing that
\begin{equation*}
  \sum_{0\le\ell\le N}\ell w^\ell
  = w\Bigl(\sum_{0\le\ell\le N}w^\ell\Bigr)'
  = \frac{w+(Nw-N-1) w^{N+1}}{(1-w)^2}.
\end{equation*}
Therefore,
\begin{equation*}
  \sum_{0\le\ell\le N}z^\ell X_\ell
  = \frac{1-wq^{-1}}{(1-w)^2} - \Biggl[
      \frac1{1-w} + \frac{N(1-q^{-1})}{1-w} + \frac{1-q^{-1}}{(1-w)^2}
    \Biggr] w^{N+1}.
\end{equation*}

If $\rho\ne0\bmod\wg^\ell$ (i.e.,\ $\ell>N$), we have
\begin{equation*}\begin{split}
  X_\ell
  &=\meas\set{(x,y)\in\oo^2:xy=\rho\bmod\wg^\ell}\\
  &=\sum_{0\le k\le N}
      \meas\set{\ord x=k ; y=\rho/x\bmod\wg^{\ell-k}}\\
  &=\sum_{0\le k\le N}q^{-k}(1-q^{-1})q^{-\ell+k}
  = (N+1)(1-q^{-1})q^{-\ell}
\end{split}\end{equation*}
and
\begin{equation*}
  \sum_{\ell>N}z^\ell X_\ell
  = (N+1)(1-q^{-1})\sum_{\ell>N}w^\ell
  = \frac{(N+1)(1-q^{-1})w^{N+1}}{1-w}.
\end{equation*}

After combining the two subsums (over $\ell\le N$ and over $\ell>N$),
the choice $\rho=t^2$ and $N=2T$ completes the proof of the isotropic case.
\end{proof}

\begin{prop}\label{p:odd:per-odd}
If $\De$ is a unit and $n=2k+2$, the local period is
\begin{multline*}
  \vol(\Th\lqu\Th K)\cdot\frac{Z(\al-n)\;L(\al-k-1,\chi)}{Z(2\al-n)}\\
  = \vol(\Th\lqu\Th K)\cdot\frac{Z\bigl((n+1)(s-1)+1\bigr)\;
          L\bigl((n+1)(s-\tfrac12)+\tfrac12,\chi\bigr)}
      {Z\bigl((n+1)(2s-1)+1\bigr)}.
\end{multline*}
\end{prop}

\begin{proof}
Let $z=q^{-\be}$, $w=zq^{-1}$, $u=z^2q^{-2}$, and $a=q^{-\al}$.  
Using $X^2$ from proposition \thref{p:odd:X2}, we obtain
\begin{equation*}
  \Pi^2(\al,\be)
  = \frac{1-\ep wq^{-1}}{(1-w)(1-\ep w)}
    \;\Biggl[\frac1{1-a}-\frac{\ep w}{1-au}\Biggr]
  = \frac{(1-\ep wq^{-1})(1+\ep aw)}{(1-w)(1-a)(1-au)}.
\end{equation*}
Therefore, the local period \eqref{e:odd:master} is
\begin{equation*}
  \frac{1+\ep aw}{Z(\al)\cdot(1-a)(1-au)}
  = \frac{Z(\al-n)\;L(\al-n+k+1,\chi)}{Z(2\al-2n+2k+2)}
\end{equation*}
and the result follows from $\al=(n+1)s$.
\end{proof}

\section{Bad odd primes}\label{s:odd:bad}

We now consider the odd primes dividing the discriminant.  We describe
their periods in terms of correction factors to those obtained
in propositions \thref{p:odd:per-even} and~\thref{p:odd:per-odd}.
In light of the discussion leading up to \eqref{e:odd:master}, we need
only consider anisotropic forms.

We assume the form $B$ is $B_0(x)+\wg B_1(y)$ for 
$(x,y)\in\oo^{n_0}\times\oo^{n_1}$, with each $B_*$ integral, regular, 
diagonal, and anisotropic.  We want to evaluate
\begin{equation*}
  X_\ell^B(t^2)
  = \meas\set{x\in\oo^{n_0},y\in\oo^{n_1}:B_0(x)+\wg B_1(y)=t^2\bmod\wg^\ell}.
\end{equation*}

If $t^2=0\bmod\wg^\ell$, the anisotropy lemma \thref{l:odd:aniso} tells us
that
\begin{equation}\label{e:odd:bad:0}\begin{split}
  X_\ell^B(0)
  &=\meas\set{x\in\oo^{n_0},y\in\oo^{n_1}:B_0(x)=\wg B_1(y)=0\bmod\wg^\ell}\\
  &=\meas\set{x\in\oo^{n_0}:B_0(x)\in\wg^\ell\oo}
    \cdot\meas\set{y\in\oo^{n_1}:B_1(y)\in\wg^{\ell-1}\oo}\\
  &= q^{-n_0\lceil\ell/2\rceil}\cdot q^{-n_1\lceil(\ell-1)/2\rceil}.
\end{split}\end{equation}

If $t^2\ne0\bmod\wg^\ell$, the anisotropy lemma, Hensel's lemma, and the
fact that $X_\ell^{B_0}(t^2)$ depends only on $T=\ord t$, tell us that
\begin{equation}\label{e:odd:bad:non0}\begin{split}
  X_\ell^B(t^2)
  &=\meas\set{x\in\oo^{n_0},y\in\oo^{n_1}:B_0(x)+\wg B_1(y)=t^2\bmod\wg^\ell}\\
  &=\meas\set{x\in\oo^{n_0},y\in\oo^{n_1}:B_0(x)=t^2\bmod\wg^\ell;
      \abs*{B_1(y)} \le \abs*{B_0(x)}}\\
  &=\meas\set{x\in\oo^{n_0}:B_0(x)=t^2\bmod\wg^\ell}\cdot q^{-n_1T}
  = X_\ell^{B_0}(t^2)\cdot q^{-n_1T}.
\end{split}\end{equation}
(If $n_0=0$, then the equation is $\wg B_1(y)=t^2\bmod\wg^\ell$, which has no
solution unless $t^2=0\bmod\wg^\ell$.)

This means we can reuse our previous results for $X_\ell^1$ and $X_\ell^2$.
Let $\De_0$ be the discriminant of $B_0$ and, from now on, $\chi$ the
quadratic character associated with $B_0$ and $\ep=\chi(\De_0)$.

\begin{prop*}
If $n_0=0$, $n_1=\ord\De=1$, and $n=2k+1$, the correction factor is
\begin{equation*}
  \frac{Z(\al-k-1)}{Z(2\al-n-1)}.
\end{equation*}
\end{prop*}

\begin{proof}
According to \eqref{e:odd:bad:0} and \eqref{e:odd:bad:non0},
\begin{align*}
  X_\ell^B(t^2)
  &=q^{-\lceil{(\ell-1)/2\rceil}},
  &\text{if $t^2=0\bmod\wg^\ell$},\\
  X_\ell^B(t^2)
  &=0,
  &\text{if $t^2\ne0\bmod\wg^\ell$}.
\end{align*}
With $z=q^{-\be}$, $w=zq^{-1}$, and $u=z^2q^{-1}$, it follows that
\begin{equation*}
  X^B(\be;t^2)
  = \sum_{0\le\ell\le 2T}z^\ell q^{-\lceil(\ell-1)/2\rceil}
  = z\sum_{-1\le\ell\le 2T-1}z^\ell q^{-\lceil\ell/2\rceil},
\end{equation*}
which, according to \eqref{e:zlqfrac}, simplifies to
\begin{equation}\label{e:bad:zlqfrac}\begin{split}
  X^B(\be;t^2)
  &=z\;\Biggl[\frac{1+w}{1-u} - \frac{w+u}{1-u}u^T + z^{-1} - u^T\Biggr]\\
  &=\frac{1+z}{1-u} - \frac{z+u}{1-u} u^T
  = \frac{1+z}{1-u}\;\Biggl[1 - \frac{u+z}{1+z} u^T\Biggr].
\end{split}\end{equation}
This fits the prototype \eqref{e:proto:X} with $v=-z$.  With $a=q^{-\al}$,
the local period \eqref{e:proto:per} is
\begin{equation*}
  \frac{Z(\al-n)}{1/(1+az)}
  = \frac{Z(\al-n) \; (1-a^2z^2)}{1-az}
  = \frac{Z(\al-n) \; Z(\al-n+k)}{Z(2\al-2n+2k)}.\qedhere
\end{equation*}
\end{proof}

\begin{prop*}
If $n_0=0$, $n_1=\ord\De=2$, and $n=2k+2$, the correction factor is
\begin{equation*}
  \frac{Z(2\al-n)\;Z(\al-k-2)}{Z(2\al-n-2)}.
\end{equation*}
\end{prop*}

\begin{proof}
According to \eqref{e:odd:bad:0} and \eqref{e:odd:bad:non0},
\begin{align*}
  X_\ell^B(t^2)
  &=q^{-2\lceil{(\ell-1)/2\rceil}},
  &\text{if $t^2=0\bmod\wg^\ell$},\\
  X_\ell^B(t^2)
  &=0,
  &\text{if $t^2\ne0\bmod\wg^\ell$}.
\end{align*}
With $z=q^{-\be}$ and $u=z^2q^{-2}$, it follows that
\begin{equation*}
  X^B(\be;t^2)
  = \sum_{0\le\ell\le 2T}z^\ell q^{-2\lceil(\ell-1)/2\rceil}.
\end{equation*}
This is the same as in the previous proposition, but with 
$q^{-2}$ instead of $q^{-1}$ (so, $u=z^2q^{-2}$ here corresponds to
$u=z^2q^{-1}$ there).  That is, the local period is
\begin{equation*}
  \frac{Z(\al-n) \; Z(\al-n+k)}{Z(2\al-2n+2k)}.\qedhere
\end{equation*}
\end{proof}

\begin{prop*}
If $n_0=n_1=\ord\De=1$ and $n=2k+2$, the correction factor is
\begin{equation*}
  \frac{Z(2\al-n)}{L(\al-k-1,\chi)},
\end{equation*}
where $\chi$ is the quadratic character associated to $B_0$.
\end{prop*}

\begin{proof}
According to \eqref{e:odd:bad:0}, \eqref{e:odd:bad:non0} and the proof of
proposition \thref{p:odd:X1},
\begin{align*}
  X_\ell^B(t^2)
  &=q^{-\ell},
  &\text{if $t^2=0\bmod\wg^\ell$},\\
  X_\ell^B(t^2)
  &=X^1_\ell(t^2)\cdot q^{-T}
  = (1+\ep) q^{-\ell},
  &\text{if $t^2\ne0\bmod\wg^\ell$}.
\end{align*}
With $z=q^{-\be}$, $w=zq^{-1}$, and $u=z^2q^{-2}=w^2$, it follows that
\begin{equation*}
  X^B(\be;t^2)
  = \sum_{0\le\ell\le 2T}z^\ell q^{-\ell}
    + (1+\ep)\sum_{\ell>2T} z^\ell q^{-\ell}
  = \frac1{1-w} + \frac{\ep w}{1-w} u^T.
\end{equation*}
With $a=q^{-\al}$, we then obtain
\begin{equation*}
  \Pi^B(\al,\be)
  = \frac1{1-w}\;\Biggl[\frac1{1-a} + \frac{\ep w}{1-au}\Biggr]
  = \frac{(1+\ep w)(1-\ep aw)}{(1-w)(1-a)(1-au)}.
\end{equation*}
That is, the local period \eqref{e:odd:master} is
\begin{equation*}\begin{split}
  \frac{1-\ep aw}{Z(\al)\cdot(1-a)(1-au)}
  = \frac{Z(\al-n)}{L(\al-n+k+1,\chi)}.\qedhere
\end{split}\end{equation*}
\end{proof}

\begin{prop*}
If $n_0=1$, $n_1=\ord\De=2$, and $n=2k+3$, the correction factor is
\begin{equation*}
  \frac1{L(\al-k-2,\chi)},
\end{equation*}
where $\chi$ is the quadratic character associated to $B_0$.
\end{prop*}

\begin{proof}
According to \eqref{e:odd:bad:0}, \eqref{e:odd:bad:non0} and the proof of
proposition \thref{p:odd:X1},
\begin{align*}
  X_\ell^B(t^2)
  &=q^{-\ell-\lceil(\ell-1)/2\rceil},
  &\text{if $t^2=0\bmod\wg^\ell$},\\
  X_\ell^B(t^2)
  &=X^1_\ell(t^2)\cdot q^{-2T}
  = (1+\ep) q^{-\ell-T},
  &\text{if $t^2\ne0\bmod\wg^\ell$}.
\end{align*}
With $z=q^{-\be}$, $w=zq^{-1}$, and $u=z^2q^{-3}=w^2q^{-1}$, it follows that
\begin{equation*}\begin{split}
  X^B(\be;t^2)
  &=\sum_{0\le\ell\le 2T}w^\ell q^{-\lceil(\ell-1)/2\rceil}
    + (1+\ep)\sum_{\ell>2T} z^\ell q^{-\ell-T}\\
  &=\frac{1+w}{1-u} - \frac{w+u}{1-u} u^T + \frac{(1+\ep)w}{1-w} u^T
  = \frac{1+w}{1-u}\;\Biggl[1 - \frac{u-\ep w}{1-\ep w} u^T\Biggr].
\end{split}\end{equation*}
(The first sum is the same as in \eqref{e:bad:zlqfrac}, with $w$ and 
$u=w^2q^{-1}$ here in place of $z$ and $u=z^2q^{-1}$ there.)  This fits the 
prototype \eqref{e:proto:X} with $v=\ep w$.  With $a=q^{-\al}$, the local 
period \eqref{e:proto:per} is
\begin{equation*}
  \frac{Z(\al-n)}{1/(1-\ep aw)}
  = \frac{Z(\al-n)}{L(\al-n+k+1,\chi)}.\qedhere
\end{equation*}
\end{proof}

\begin{prop*}
If $n_0=2$, $n_1=\ord\De=1$, and $n=2k+3$, the correction factor is
\begin{equation*}
  \frac1{Z(\al-k-2)}.
\end{equation*}
\end{prop*}

\begin{proof}
According to \eqref{e:odd:bad:0}, \eqref{e:odd:bad:non0} and the proof of
proposition \thref{p:odd:X2} (anisotropic case),
\begin{align*}
  X_\ell^B(t^2)
  &=q^{-\ell-\lceil\ell/2\rceil},
  &\text{if $t^2=0\bmod\wg^\ell$},\\
  X_\ell^B(t^2)
  &=X^1_\ell(t^2)\cdot q^{-T}
  = (1+q^{-1}) q^{-\ell-T},
  &\text{if $t^2\ne0\bmod\wg^\ell$}.
\end{align*}
With $z=q^{-\be}$, $w=zq^{-1}$, and $u=z^2q^{-3}=w^2q^{-1}$, we have
\begin{equation*}
  \sum_{0\le\ell\le2T} z^\ell X_\ell^B(t^2)
  = \sum_{0\le\ell\le2T} w^\ell q^{-\lceil\ell/2\rceil}
  = \frac{1+wq^{-1}}{1-u} - \frac{wq^{-1}+u}{1-u} u^T,
\end{equation*}
according to \eqref{e:zlqfrac} (with $w$, $wq^{-1}$, and $u=w^2q^{-1}$ 
here in place of $z$, $w=zq^{-1}$, and $u=z^2q^{-1}$ there).  We have also
\begin{equation*}
  \sum_{\ell>2T} z^\ell X_\ell^B(t^2)
  = (1+q^{-1})q^{-T} \sum_{\ell>2T} w^\ell
  = \frac{(1+q^{-1}) w}{1-w} u^T,
\end{equation*}
leading to
\begin{equation*}\begin{split}
  X^B(\be;t^2)
  &=\frac{1+wq^{-1}}{1-u} - \frac{wq^{-1}+u}{1-u}u^T
    + \frac{w+wq^{-1}}{1-w} u^T\\
  &=\frac{1+wq^{-1}}{1-u}\;\Biggl[1 - \frac{u-w}{1-w} u^T\Biggr].
\end{split}\end{equation*}
This fits the prototype \eqref{e:proto:X} with $v=w$.  With $a=q^{-\al}$, the
local period \eqref{e:proto:per} is
\begin{equation*}
  \frac{Z(\al-n)}{1/(1-aw)}
  = \frac{Z(\al-n)}{Z(\al-n+k+1)}.\qedhere
\end{equation*}
\end{proof}

\begin{prop*}
If $n_0=n_1=\ord\De=2$ and $n=2k+4$, the correction factor is
\begin{equation*}
  \frac{Z(2\al-n)}{Z(\al-k-3)}.
\end{equation*}
\end{prop*}

\begin{proof}
According to \eqref{e:odd:bad:0}, \eqref{e:odd:bad:non0} and the proof of
proposition \thref{p:odd:X2} (anisotropic case),
\begin{align*}
  X_\ell^B(t^2)
  &=q^{-2\ell},
  &\text{if $t^2=0\bmod\wg^\ell$},\\
  X_\ell^B(t^2)
  &=X^1_\ell(t^2)\cdot q^{-2T}
  = (1+q^{-1})  q^{-\ell-2T},
  &\text{if $t^2\ne0\bmod\wg^\ell$}.
\end{align*}
With $z=q^{-\be}$, $w=zq^{-1}$, and $u=z^2q^{-4}=w^2q^{-2}$, we have
\begin{equation*}\begin{split}
  X^B(\be;t^2)
  &=\sum_{0\le\ell\le2T}z^\ell q^{-2\ell}
    +(1+q^{-1})\sum_{\ell>2T}z^\ell q^{-\ell-2T}\\
  &=\frac{1-wq^{-1}u^T}{1-wq^{-1}} + \frac{w+wq^{-1}}{1-w} u^T
  = \frac1{1-wq^{-1}}\;\Biggl[1 - \frac{u-w}{1-w} u^T\Biggr].
\end{split}\end{equation*}
This fits the prototype \eqref{e:proto:X} with $v=w$.  With $a=q^{-\al}$, the
local period \eqref{e:proto:per} is
\begin{equation*}
  \frac{Z(\al-n)}{1/(1-aw)}
  = \frac{Z(\al-n)}{Z(\al-n+k+1)}.\qedhere
\end{equation*}
\end{proof}

\section{The general case}\label{s:gen}

We resume the discussion and the notation from section~\ref{s:setup}.  
We saw that
\begin{equation*}\tag{\ref{e:auto-H}}
  (E_\ph,F)_H
  = (\eta,F)_\Th \cdot \int_{\Sa\lqu\Ha}\ph_\om \cdot f,
\end{equation*}
where $f$ is a spherical vector in $\Ind_\Th^H1$ normalized by $f(1)=1$.
We want to express this global integral as a product of local integrals.

If $F$ generates an irreducible representation, so does $f$ generate an 
irreducible representation $\pi=\otimes_v\pi_v$.  If the period is nonzero, 
then each $\pi_v$ is a quotient of a degenerate unramified principal series 
with respect to the Levi component $M^Q$ and with character 
$m_\la\mapsto\abs\la^{\be_v}$.  Let $f_v$ be a generator of $\pi_v$ 
normalized by $f_v(1)=1$.  The intertwining between models of $\pi$ takes 
$f$ to a multiple of $\otimes_vf_v$, so
\begin{equation*}
  \int_{\Sa\lqu\Ha}\ph_\om\cdot f
  = C_f \cdot \prod_v \int_{\Sv\lqu\Hv}\ph_{\om,v}\cdot f_v
\end{equation*}
for some constant $C_f$.

We now work locally and suppress $v$.  Let $\psi_\be$ be the principal
series; in particular, $\psi_\be(m_\la)=\abs\la^\be$.  The function $f$ 
generates a model of $\pi$ consisting of left $\Theta$-- and right 
$K$--invariant functions.  The restriction $\Res_Q^H\pi$ may be 
modelled by the very same functions.  Thus, it is irreducible too and 
isomorphic to the principal series representation.  But we may model both 
of them with left $\Th\cap Q$-- and right $K$--invariant functions and 
$f(1)=\psi_\be(1)$; therefore, $\psi_\be$ is the restriction of $f$ to $Q$.

\begin{prop}\label{p:general}
The local integral at odd primes is given by
\begin{equation*}
  \int_{\Th\lqu H}\ph_\om\cdot f
  = \frac{\vol(\Th\lqu\Th K)\cdot\Pi^n(\al-\be-n,\be)}{Z(\al)\cdot X^n(0;1)}.
\end{equation*}
\end{prop}

\begin{proof}
This is the same computation we did in the proof of 
proposition \thref{p:odd:lpi}; by the exact same reasoning and with the same
notation as there, the local integral is
\begin{equation*}
  \frac{\vol(\Th\lqu\Th K)}{Z(\al)\cdot X^n(0;1)}
      \int_{(\Th\cap Q)\lqu Q}\ph_\om\cdot f
  = \frac{\vol(\Th\lqu\Th K)}{Z(\al)\cdot X^n(0;1)}
      \int_{(\Th\cap Q)\lqu Q}\ph_\om\cdot\abs\la^\be.
\end{equation*}
But that last integral is the very same one computed from 
\eqref{e:odd:triple} onward.
\end{proof}

\appendix
\section{The constant term}\label{s:const}

The functional equation of an Eisenstein series is inherited by its periods.
As the equation involves the constant term, we describe briefly how that
term may be obtained using the methods we have used so far.

With $\vol(\Nq)=1$, the 
constant term of $E_s$ is
\begin{equation*}
  c\,E_s(m)
  = \phs(m) + \int_{\Na}\phs(\xi nm)\dd n
  = \abs\la^\al + \int_{\Na}\phs(\xi nm)\dd n,
\end{equation*}
where $\xi$ is the long Weyl element in $P\lqu G/P$.
The second summand factors over primes.

At a place $v$ (suppressed from this point onward), 
express the form on $k^{n+3}=(k\cdot e')\oplus k^{n+1}\oplus(k\cdot e)$ by
\begin{equation*}
  \begin{pmatrix}&&1\\&B&\\1&&\end{pmatrix}
\end{equation*}
and write
\begin{equation*}
  n = \begin{pmatrix}1&a&-\frac12B(a)\\&\id&*\\
    \phantom{\frac12B(a)}&&1\end{pmatrix}
  \qquad\text{and}\qquad
  m = \begin{pmatrix}\la&&\\&\theta&\\&&1/\la\end{pmatrix},
\end{equation*}
with $\theta\in\Th=\gO(B)$.  We have
\begin{equation*}\begin{split}
  \int_N \ph(\xi nm)\dd n
  &=\int_N\int_\kx\abs t^\al\;\Phi(te\cdot\xi nm)\dd t\dd n
  = \int_N\int_\kx\abs t^\al\;\Phi(te'\cdot nm)\dd t\dd n\\
  &=\int_N\int_\kx\abs t^\al\;\Phi(te'\cdot mn)\dd t\cdot\,\de_P(m)\dd n\\
  &=\int_N\int_\kx\abs t^\al\;\abs\la^{n+1}\;\Phi(\la te'\cdot n)\dd t\dd n.
\end{split}\end{equation*}
At this point, we disregard the normalization of the integral, and proceed.
\begin{equation*}\begin{split}
  \int_N \ph(\xi nm)\dd n
  &=\abs\la^{(n+1)-\al}\cdot\int_\kx\abs t^\al\int_N
      \Phi(te'\cdot n)\dd n\dd t\\
  &=\abs\la^{(n+1)-\al}\cdot\int_\kx\abs t^\al\int_{k^{n+1}}
      \Phi\Bigl(t,at,-\tfrac1{2t}B(at)\Bigr)\dd a\dd t\\
  &=\abs\la^{(n+1)-\al}\int_\kx\abs t^{\al-(n+1)}\int_{k^{n+1}}
    \Phi\Bigl(t,a,-\tfrac1{2t}B(a)\Bigr)\dd a\dd t.
\end{split}\end{equation*}

Using the same notation we used before
\begin{equation*}
  \int_N\ph(\xi nm)\dd n
  = \abs\la^{(n+1)(1-s)}\cdot X^{n+1}\bigl(\al-(n+1);0\bigr).
\end{equation*}
Therefore, recalling $\al=(n+1)s$, 
\begin{equation*}
  \int_N \phs(\xi nm)\dd n
  = \abs\la^{(n+1)(1-s)}\cdot
    \frac{X^{n+1}(-(n+1)(1-s);0)}{Z\bigl((n+1)s\bigr)}.
\end{equation*}

With the conventions of the strategy overview at the end of 
section~\ref{s:odd},
\begin{equation*}
  X(\be;0)
  = \lim_{a\to1}\Bigl((1-a)\sum_{T\ge0}a^T X(\be;\wg^{2T})\Bigr)
  = -\Res_{a=1}R(z,w,u,a).
\end{equation*}

\nocite{Ca,GiJiSo,GiPSRa,GiRaSo,Gr,
  Ig,IwSa,JaSh1,JaSh2,JiSo1,JiSo2,Ki,Om,PlRa,PS,Sc,SuZh}


\begin{bibdiv}


\begin{biblist}

\bib{AiGoRaSc}{article}{
  label={AGRS08},
  author={Aizenbud, Avraham},
  author={Gourevitch, Dmitry},
  author={Rallis, Steve},
  author={Schiffmann, G{\'e}rard},
  title={Multiplicity one theorems},
  journal={Ann. of Math. (2)},
  volume={172},
  date={2010},
  number={2},
  pages={1407--1434},
  review={\MR{2680495 (2011g:22024)}},
}

\bib{AiGoSa}{article}{
  author={Aizenbud, Avraham},
  author={Gourevitch, Dmitry},
  author={Sayag, Eitan},
  title={{$(O(V\oplus F),O(V))$} is a {G}elfand pair for any quadratic space {$V$} over a local field {$F$}},
  date={2009},
  journal={Math. Z.},
  volume={261},
  number={2},
  pages={239\ndash 244},
  review={\MR{2457297 (2010a:22020)}},
}

\bib{Ca}{book}{
  author={Cassels, J. W.~S.},
  title={Rational quadratic forms},
  series={London Mathematical Society Monographs},
  publisher={Academic Press Inc.},
  address={London},
  date={1978},
  volume={13},
  review={\MR{522835 (80m:10019)}},
}

\bib{Ga}{article}{
  author={Garrett, Paul},
  title={Euler factorization of global integrals},
  pages={35\ndash 101},
  book={
    title={Automorphic forms, automorphic representations, and arithmetic},
    editor={Doran, Robert~S.},
    editor={Dou, Ze-Li},
    editor={Gilbert, George~T.},
    series={Proc. Sympos. Pure Math.},
    publisher={Amer. Math. Soc.},
    address={Providence, R.I.},
    volume={66.2}, },
  date={1999},
  conference={
    title={NSF\ndash CBMS Regional Conference in Mathematics on Euler Products and Eisenstein Series},
    address={Texas Christian University, Fort Worth, Tex.},
    date={May 20\ndash 24, 1996}, },
  review={\MR{1703758 (2000m:11043)}},
}

\bib{GePSRa}{book}{
   author={Gelbart, Stephen},
   author={Piatetski-Shapiro, Ilya},
   author={Rallis, Stephen},
   title={Explicit constructions of automorphic $L$-functions},
   series={Lecture Notes in Mathematics},
   volume={1254},
   publisher={Springer-Verlag},
   place={Berlin},
   date={1987},
   review={\MR{892097 (89k:11038)}},
}

\bib{GiJiSo}{article}{
  author={Ginzburg, David},
  author={Jiang, Dihua},
  author={Soudry, David},
  title={Poles of {$L$}-functions and theta liftings for orthogonal groups},
  date={2009},
  journal={J. Inst. Math. Jussieu},
  volume={8},
  number={4},
  pages={693\ndash741},
  review={\MR{2540878}},
}

\bib{GiPSRa}{article}{
  author={Ginzburg, D.},
  author={Piatetski-Shapiro, I.},
  author={Rallis, S.},
  title={{$L$} functions for the orthogonal group},
  date={1997},
  journal={Mem. Amer. Math. Soc.},
  volume={128},
  number={611},
  review={\MR{1357823 (98m:11041)}},
}

\bib{GiRaSo}{article}{
  author={Ginzburg, David},
  author={Rallis, Stephen},
  author={Soudry, David},
  title={Periods, poles of {$L$}-functions and symplectic-orthogonal theta lifts},
  date={1997},
  journal={J. Reine Angew. Math.},
  volume={487},
  pages={85\ndash 114},
  review={\MR{1454260 (98f:11046)}},
}

\bib{Gr}{article}{
  author={Gross, Benedict~H.},
  title={Some applications of {G}elfand pairs to number theory},
  date={1991},
  journal={Bull. Amer. Math. Soc. (N.S.)},
  volume={24},
  number={2},
  pages={277\ndash 301},
  review={\MR{1074028 (91i:11055)}},
}

\bib{GrPr1}{article}{
  author={Gross, Benedict~H.},
  author={Prasad, Dipendra},
  title={Test vectors for linear forms},
  date={1991},
  journal={Math. Ann.},
  volume={291},
  number={2},
  pages={343\ndash 355},
  review={\MR{1129372 (92k:22028)}},
}

\bib{GrPr2}{article}{
  author={Gross, Benedict~H.},
  author={Prasad, Dipendra},
  title={On the decomposition of a representation of {${\rm SO}\sb n$} when restricted to {${\rm SO}\sb {n-1}$}},
  date={1992},
  journal={Canad. J. Math.},
  volume={44},
  number={5},
  pages={974\ndash 1002},
  review={\MR{1186476 (93j:22031)}},
}

\bib{GrPr3}{article}{
  author={Gross, Benedict~H.},
  author={Prasad, Dipendra},
  title={On irreducible representations of {${\rm SO}\sb {2n+1} \times {\rm SO}\sb {2m}$}},
  date={1994},
  journal={Canad. J. Math.},
  volume={46},
  number={5},
  pages={930\ndash 950},
  review={\MR{1295124 (96c:22028)}},
}

\bib{GrRe}{article}{
  author={Gross, Benedict~H.},
  author={Reeder, Mark},
  title={From {L}aplace to {L}anglands via representations of orthogonal groups},
  date={2006},
  journal={Bull. Amer. Math. Soc. (N.S.)},
  volume={43},
  number={2},
  pages={163\ndash 205},
  review={\MR{2216109 (2007a:11159)}},
}

\bib{IcIk}{article}{
  author={Ichino, Atsushi},
  author={Ikeda, Tamotsu},
  title={On the periods of automorphic forms on special orthogonal groups and the {G}ross--{P}rasad conjecture},
  date={2010},
  journal={Geom. Funct. Anal.},
  volume={19},
  number={5},
  pages={1378\ndash1425},
  review={\MR{2585578}},
}

\bib{Ig}{book}{
  author={Igusa, Jun-ichi},
  title={An introduction to the theory of local zeta functions},
  series={AMS/IP Studies in Advanced Mathematics},
  publisher={Amer. Math. Soc.},
  address={Providence, R.I.},
  date={2000},
  volume={14},
  review={\MR{1743467 (2001j:11112)}},
}

\bib{IwSa}{article}{
  author={Iwaniec, H.},
  author={Sarnak, P.},
  title={Perspectives on the analytic theory of {$L$}-functions},
  date={2000},
  book={
    title={GAFA 2000},
    editor={Alon, N.},
    editor={Bourgain, J.},
    editor={Connes, A.},
    editor={Gromov, M.},
    editor={Milman, V.},
    publisher={Birkh\"auser},
    address={Basel}, },
  conference={
    title={Visions in Mathematics, towards 2000},
    address={Tel Aviv University},
    date={August 25\ndash September 3, 1999}, },
  pages={705\ndash 741},
  note={Geom. Funct. Anal., Special Volume, Part II},
  review={\MR{1826269 (2002b:11117)}},
}

\bib{JaLaRo}{article}{
  author={Jacquet, Herv{\'e}},
  author={Lapid, Erez},
  author={Rogawski, Jonathan},
  title={Periods of automorphic forms},
  date={1999},
  journal={J. Amer. Math. Soc.},
  volume={12},
  number={1},
  pages={173\ndash 240},
  review={\MR{1625060 (99c:11056)}},
}

\bib{JaSh1}{article}{
  author={Jacquet, H.},
  author={Shalika, J.~A.},
  title={On {E}uler products and the classification of automorphic representations. {I}},
  date={1981},
  journal={Amer. J. Math.},
  volume={103},
  number={3},
  pages={499\ndash 558},
  review={\MR{618323 (82m:10050a)}},
}

\bib{JaSh2}{article}{
  author={Jacquet, H.},
  author={Shalika, J.~A.},
  title={On {E}uler products and the classification of automorphic forms. {II}},
  date={1981},
  journal={Amer. J. Math.},
  volume={103},
  number={4},
  pages={777\ndash 815},
  review={\MR{623137 (82m:10050b)}},
}

\bib{Ji}{article}{
  author={Jiang, Dihua},
  title={Periods of automorphic forms},
  date={2007},
  book={
    title={Proceedings of the International Conference on Complex Geometry and Related Fields},
    editor={Yau, Stephen S.-T.},
    editor={Chen, Zhijie},
    editor={Wang, Jianpan},
    editor={Ten, Sheng-Li},
    series={AMS/IP Stud. Adv. Math.},
    volume={39},
    publisher={Amer. Math. Soc.},
    address={Providence, R.I.}, },
  pages={125\ndash 148},
  review={\MR{2338623 (2008k:11052)}},
}

\bib{JiSo1}{article}{
  author={Jiang, Dihua},
  author={Soudry, David},
  title={The local converse theorem for {${\rm SO}(2n+1)$} and applications},
  date={2003},
  journal={Ann. of Math. (2)},
  volume={157},
  number={3},
  pages={743\ndash 806},
  review={\MR{1983781 (2005b:11193)}},
}

\bib{JiSo2}{article}{
  author={Jiang, Dihua},
  author={Soudry, David},
  title={Generic representations and local {L}anglands reciprocity law for {$p$}-adic {${\rm SO}\sb {2n+1}$}},
  date={2004},
  book={
    title={Contributions to automorphic forms, geometry, and number theory},
    editor={Hida, Haruzo},
    editor={Ramakrishnan, Dinakar},
    editor={Shahidi, Freydoon},
    publisher={Johns Hopkins Univ. Press},
    address={Baltimore, Md.}, },
  pages={457\ndash 519},
  review={\MR{2058617 (2005f:11272)}},
}

\bib{KaMuSu}{article}{
   author={Kato, Shin-ichi},
   author={Murase, Atsushi},
   author={Sugano, Takashi},
   title={Whittaker--Shintani functions for orthogonal groups},
   journal={Tohoku Math. J. (2)},
   volume={55},
   date={2003},
   number={1},
   pages={1--64},
   review={\MR{1956080 (2003m:22020)}},
}

\bib{Ki}{article}{
  author={Kim, Henry~H.},
  title={On local {$L$}-functions and normalized intertwining operators},
  date={2005},
  journal={Canad. J. Math.},
  volume={57},
  number={3},
  pages={535\ndash 597},
  review={\MR{2134402 (2006a:11063)}},
}

\bib{LaOf}{article}{
  author={Lapid, Erez},
  author={Offen, Omer},
  title={Compact unitary periods},
  date={2007},
  journal={Compos. Math.},
  volume={143},
  number={2},
  pages={323\ndash 338},
  review={\MR{2309989 (2008g:11091)}},
}

\bib{LaRo}{article}{
  author={Lapid, Erez},
  author={Rogawski, Jonathan},
  title={Periods of {E}isenstein series},
  date={2001},
  journal={C. R. Acad. Sci. Paris S\'er. I Math.},
  volume={333},
  number={6},
  pages={513\ndash 516},
  review={\MR{1860921 (2002k:11072)}},
}

\bib{MuSu}{article}{
   author={Murase, Atsushi},
   author={Sugano, Takashi},
   title={Shintani function and its application to automorphic $L$-functions
   for classical groups. I. The case of orthogonal groups},
   journal={Math. Ann.},
   volume={299},
   date={1994},
   number={1},
   pages={17--56},
   review={\MR{1273075 (96c:11054)}},
}

\bib{Om}{book}{
  author={O'Meara, O.~Timothy},
  title={Introduction to quadratic forms},
  series={Classics in Mathematics},
  publisher={Springer-Verlag},
  address={Berlin},
  date={2000},
  note={Reprint of the 1973 edition},
  review={\MR{1754311 (2000m:11032)}},
}

\bib{PS}{article}{
  author={Pjateckij-{\v {S}}apiro, I.~I.},
  title={Euler subgroups},
  pages={597\ndash 620},
  book={
    editor={Gelfand, I.~M.},
    title={Lie groups and their representations},
    publisher={John Wiley \& Sons},
    address={New York}, },
  date={1975},
  conference={
    title={Summer School on Group Representations of the Bolyai J{\'a}nos Mathematical Society},
    address={Budapest},
    date={August 16\ndash September 3, 1971}, },
  review={\MR{0406935 (53 \#10720)}},
}

\bib{PlRa}{book}{
  author={Platonov, Vladimir},
  author={Rapinchuk, Andrei},
  title={Algebraic groups and number theory},
  series={Pure and Applied Mathematics},
  publisher={Academic Press Inc.},
  address={Boston, Mass.},
  date={1994},
  volume={139},
  note={Translated from the 1991 Russian original by Rachel Rowen},
  review={\MR{1278263 (95b:11039)}},
}

\bib{Ra}{book}{
   author={Rallis, Stephen},
   title={$L$-functions and the oscillator representation},
   series={Lecture Notes in Mathematics},
   volume={1245},
   publisher={Springer-Verlag},
   place={Berlin},
   date={1987},
   review={\MR{887329 (89b:11046)}},
}

\bib{Sc}{book}{
  author={Scharlau, Winfried},
  title={Quadratic and {H}ermitian forms},
  series={Grundlehren der Mathematischen Wissenschaften},
  publisher={Springer-Verlag},
  address={Berlin},
  date={1985},
  volume={270},
  review={\MR{770063 (86k:11022)}},
}

\bib{SuZh}{article}{
  author={Sun, Binyong},
  author={Zhu, Chen-Bo},
  title={Multiplicity one theorems: the archimedean case},
  journal={Ann. of Math. (2)},
  status={to appear},
  eprint={ arXiv:0903.1413v2 [math.RT] },
}

\bib{Wei1}{book}{
  author={Weil, Andr{\'e}},
  title={Adeles and algebraic groups},
  series={Progress in Mathematics},
  publisher={Birkh\"auser},
  address={Boston, Mass.},
  date={1982},
  volume={23},
  review={\MR{670072 (83m:10032)}},
}


\end{biblist}

\end{bibdiv}
\end{document}